\documentclass[final]{elsart}
\usepackage{graphicx,amssymb,amsmath,amsfonts}

\newcommand{\De}{\Delta}
\newcommand{\de}{\partial}
\newcommand{\lu}[1]{\;{}^#1\!}
\def\norm#1{\Vert #1 \Vert_1}
\def\mod#1{\left| #1 \right|}
\def\normtot#1{\Vert #1 \Vert_{\infty,1}}
\newcommand{\F}{\mathcal{F}}

\begin{document}
\begin{frontmatter}

\title{A finite difference method for \\ piecewise deterministic Markov processes}
\author{Mario Annunziato}
\ead{email:mannunzi@unisa.it}
\address{Dipartimento di Matematica e Informatica \\
Universit\`a degli Studi di Salerno\\ 
Via Ponte Don Melillo 84084 (SA)}

\begin{abstract}
An extension of non-deterministic processes driven by the random
telegraph signal is introduced in the framework of \textit{piecewise deterministic Markov processes} \cite{dav},
including a broader category of random systems.
The corresponding Liouville-Master Equation is established and the upwind method is applied to numerical 
calculation of the distribution function.
The convergence of the numerical solution is proved under an appropriate Courant-Friedrichs-Lewy condition. 
The same condition  preserve the non-decreasing property of the calculated
distribution function. Some numerical tests are presented.
\end{abstract}%

\end{frontmatter}

\section{Introduction}

Modeling non-deterministic physical systems is a modern subject of 
research that covers many different disciplines. For models that
use an \textit{a priori} source of randomness, the Gaussian white
noise is the most common ingredient. Nevertheless there is a growing interest 
in searching alternatives to the white noise. There are two main characteristics,
often considered unrealistic for modeling, that motivate these researches:
the white noise induces (i) infinitely many randomness for each small interval 
of time and (ii) unlimited fast fluctuations, although with an exponentially decreasing probability.
The random telegraph signal (named also binary or dichotomous process) is 
the most elementary assumption that does not suffer of these
pathologies: it is a signal switching randomly between two values,
with a determined law for switching times.
As examples of application to science and engineering,
we mention: anomalous diffusion \cite{bol:gri:wes}, reaction-diffusion \cite{hor}, scattering of radiation  \cite{jak:ren}, biological dispersal \cite{oth:dun:alt}, for diffusive processes, 
and non-Maxwellian equilibrium \cite{mario,cac,fil:hon,kit,bro:han},
diagnostic technique for semiconductor lasers \cite{jak:rid}, 
filtered telegraph signal 
\footnote{The author acknowledge K.D. Jakeman and E. Ridley for having
advised of Pawula's works.} \cite{paw:ric,jak:rid},  harmonic
oscillators \cite{mas:por}, for systems having an equilibrium.

For an introduction, let us consider a dissipative process subject to a
noised input:
\begin{equation}
\frac{dX}{dt} = - X(t) +  \xi(t)
\label{FRTP}
\end{equation}
where $X(t)$ is a dynamical variable and $\xi(t)$ represents the noise. 
When $\xi(t)$ is the Gaussian white noise, $X(t)$ is an 
Ornstein-Uhlenbeck process, and the associated equation
$\de_t p(x,t) = \de_x (x\,p(x,t)) +(1/2)\de_x^2 p(x,t)$,
describing the density distribution function $p(x,t)$, results by
the well known Fokker-Planck equation \cite{gar}.
When $\xi(t)$ is taken as the random telegraph signal, the same equation
describes a filtered random telegraph process 
\cite{paw:ric,jak:ren} or, equivalently, a Langevin equation \cite{mario,cac} subject to a
dichotomous noise.
$\xi(t)$ alternately takes on values $\pm 1$, 
with an exponential (or Poisson) switch 
probability density function of the form:
 $\mu \, e^{-\mu \tau}$, where $\mu^{-1}$ is the
average time $\langle\tau\rangle$ between switches.
With this assumption,
we see that the Eq. (\ref{FRTP}) can be integrated as an ordinary 
differential equation, provided that not any switching event 
happen inside the integration interval. Therefore, with the exceptions of switching times, 
the process is deterministic, composed of pieces of increasing and decreasing exponentials.
Anyway, the whole resulting process $X(t)$ is not
deterministic, it represents a random sample path in a probability
space, and  the statistical properties can be described 
using the characteristic functional method \cite{mor,cac,bud:cac,jak:rid}
or the associated probability densities distribution $p^{\pm}(x,t)$,
governed by a Liouville-Master Equation
\footnote{In general, equations for density probability of random processes
are derived from the Chapman-Kolmogorov equation. As discussed in Ref. \cite{gar}
the same equation becomes a Liouville equation, in absence of randomness, and
a Master Equation, if only jump processes are involved. We use both terms
in order to stress the deterministic and the random character of the processes
considered here.} 
of the following form \cite{won,jak:rid,mario}:
\begin{equation}
\left\{
\begin{array}{l}
\de_t p^+ -(x-1) \de_x p^+ = (1-\mu) p^+ + \mu p^-  \\
\de_t p^- -(x+1) \de_x p^- = \mu p^+ + (1-\mu) p^-.
\end{array}
\right.
\label{ME_Filter_RTS}
\end{equation}

%

This kind of process is a case of a more general class of \textit{piecewise-deterministic processes} that 
were formalized by Davis \cite{DAV,dav}. 
Following this framework, in Sect. \ref{sec:extend} we define an extension 
of Eq. (\ref{FRTP}), generating continuous piecewise deterministic processes, 
and establish the corresponding  Liouville-Master Equation, like Eq. (\ref{ME_Filter_RTS}).
Sect. \ref{sec:numsol} is devoted to investigate the properties of the 
upwind method used to solve the Liouville-Master Equation. 
In Sect. \ref{numtest} the numerical solution of a dissipative model is plotted and 
compared to Monte Carlo's histograms of the same process.

\section{Continuous piecewise-deterministic Markov processes}
\label{sec:extend}

The process described by Eq. (\ref{FRTP}), rather than a single equation driven by a binary noise, 
can be regarded as composed of pieces of
two deterministic equations having two velocity states. 
It is clear that following this point of view, we can perform some extensions:
(i) on the number of states of the system, (ii) on the laws  governing 
the velocity of the dynamical variable $X(t)$, (iii) on the
statistics of the generating randomness.
Based on these ideas we can formalize a piecewise deterministic process as follows:
\begin{enumerate}
\item[(a)]  The dynamics is described by the equation:
\begin{equation}
\frac{dX}{dt} = A_{s}(X)
\label{RTE_estesa_1}
\end{equation}
where $A_{s}$ is a function chosen randomly in a set of $\lbrace A_1,\ldots,A_S\rbrace$ functions. 
Given $A_s$, we say that the dynamics is in the state $s$.
We require that $A_s(x)$ be Lipschitz continuous, so that 
for fixed $s$, $X(t)$ exists and is unique. 
\item[(b)] The initial condition is settled either by the Cauchy problem to Eq. (\ref{RTE_estesa_1}), 
i.e. $X(t_0)=X_0$, and by the initial state of the process. 
%
\item[(c)] States are characterised by a Poisson statistics of 
\textit{switching times}: 
\begin{equation}
\mu_s \, e^{-\mu_s t}
\label{RTE_estesa_2}
\end{equation}
with $\mu_s \in \{\mu_1,\ldots,\mu_S\}$. 
\item[(d)]  At every \textit{switching time} the dynamics can change from the 
state $j$ ($A_j$) to the state $i$ ($A_i$) according to a 
transition probability matrix (or stochastic matrix):
\begin{equation}
q_{ij}
\label{RTE_estesa_3}
\end{equation}
having the following fundamental \cite{gar} properties:  $0 \leq q_{ij}\leq 1$ and $\sum_{i=1}^S q_{ij} =1$, for
each $i,j \in 1,\ldots,S$.
\end{enumerate}
We can easily visualize this process as a sequence of solutions of 
Eq. (\ref{RTE_estesa_1}), where the end point of each is just 
the initial condition for the successive integral curve, so that
 this process is continuous, because at switching times
only the first derivative of $X(t)$ is involved. Assumption (c) makes 
the process be Markovian \cite{DAV,dav}, because the exponential statistics, regarded as
generated by an infinitesimal time process, corresponds to a constant switching rate.

\subsection{The Liouville-Master Equation}

The statistical description of the process (a,b,c,d) 
is performed by the Liouville-Master Equation for
the distribution functions $F_s(x,t)=\int_{-\infty}^x p_s(x',t) dx' $, 
for each state $s=\left\lbrace 1, \ldots, S \right\rbrace $ of the system:
\begin{equation}
      \de_t F_s(x,t) + A_s(x)\, \de_x F_s(x,t)  =
      \sum_{j=1}^{S} Q_{sj}  F_j(x,t) ,
\label{ME}
\end{equation}
where $Q_{sj}=(q_{sj}-\delta_{sj})\mu_s$, and $\delta_{sj}$ is the Kronecker's symbol. 
Eq. (\ref{ME}) is written for $F_s(x,t)$ rather than the densities $p_s(x,t)$, 
but they are connected by a derivative. 
We remind that $p_s(x,t) \geq 0$ and $F_s(x,t)$ is non-decreasing in $x$.
The Eq. (\ref{ME}) is solved for the Cauchy initial conditions:
\begin{equation}
F_s(x,0) = F_{0s}(x) = \int_{-\infty}^x p_{0s}(x') dx'
\label{cauchy}
\end{equation}
given for each state $s$ of the system.
Our aim is to find solutions $F_1(x,t),\cdots,$ $F_S(x,t)$ of 
Eq. (\ref{ME}), and then the total distribution function: 
$\F(x,t) = \sum_{s=1}^{S} F_s(x,t)$, representing  
the distribution of the process at time $t$, regardless of the state of the equation. 
Note that the conservation of probability measure requires that both 
conditions $lim_{x\rightarrow \infty} \F(x,t) = 1$
and $lim_{x\rightarrow -\infty} \F(x,t) = 0$, must be satisfied.
This means that Eq. (\ref{ME}) is also subject to a limit problem.
We note that if $A_1(x) = -x+1$, $A_2(x) = -x-1$, $\mu_1=\mu_2=\mu$,
$q_{ij}=1$ for $i\neq j$ and $q_{ij}=0$ for $i=j$, then
the filtered random telegraph process (\ref{FRTP}) and its
associated Liouville-Master Equation (\ref{ME_Filter_RTS}) are recovered.

This system of equation is linear hyperbolic, with non constants coefficients. 
In what follows we focus our attention on solutions $\F(x,t)$
having an equilibrium, but the numerical scheme can be extended 
to diffusion processes too.
Conditions for the existence of equilibrium solution can be conjectured
by using simple dynamical considerations \cite{mor}.
If all dynamical equations (\ref{RTE_estesa_1}) own only attraction points and
all of these are contained into the intersection of the basin of 
attraction of each, then a process starting from this region
will never escape. Whence, there should exists a region $\Omega$ where the 
process is confined and a stationary distribution $\F_{eq}(x)=\lim_{t\rightarrow\infty}\F(x,t)$ 
exists.

\section{The numerical solution}
\label{sec:numsol}

Due to his hyperbolic nature, we use the \textit{upwind} scheme 
\cite{mor:may,ran} to numerically 
solve Eq. (\ref{ME}). 
We define the uniform discrete mesh points ($x_k,t_i$) with steps $\De x$ and $\De t$ on the domain $(\Omega \times [0,T])$, so that the coefficient 
$A_l(x)$ is written as $ \lu lA_{k} = A_l(x_k)$ and the discrete solution as
$ \lu l{F}_k^{i} = F_l(x_k,t_i) $.
The upwind scheme takes the discrete derivative $x$-direction according 
to the sign of $A_s(x)$, so that we obtain the following difference discrete 
scheme for $\lu l{F}_k^{i}$:
\begin{equation}
\begin{array}{c}
\lu lF_k^{i+1} = \lu l{F}_k^{i} + \Delta t \left(-\lu lA_{k} 
             \frac{\lu lF_{k+\nu}^{i} - \lu lF_{k+\nu-1}^{i}}{\Delta x}
		+ \sum_s Q_{ls}  \; \lu s{F}_k^{i} \right) \\
\mbox{where $\nu =$}
\left\lbrace 
\begin{array}{l}
1 \;\;\;\mbox{ if } \;\;\; \lu lA_{k} \leq 0 \\
0 \;\;\;\mbox{ if } \;\;\; \lu lA_{k} \geq 0 \\
\end{array}
\right.
\end{array}
\label{upwind}
\end{equation} 
$\nu=1$  corresponds to take in Eq. (\ref{upwind}) 
the spatial discrete right-derivative; $\nu=0$ the left one.
This is an explicit scheme of first order, where the solutions $\lu lF_k^n$, at time
step $n$, are calculated from the 
Cauchy initial conditions $\lu lF_k^0 = F_{0l}(x_k)$.

\subsection{Convergence analysis}

Stability conditions for the upwind method are known for
linear hyperbolic systems as Courant-Friedrichs-Lewy (CFL).
If $Q=0$ the numerical scheme (\ref{upwind}) is stable when $\mid A_s \;\De t/\De x \mid\leq 1$
is fulfilled. If $Q\neq 0$ we can use the approach of Ref. \cite{mor:may}
to state that the discrete operator of (\ref{upwind}) 
is limited by an exponential grow in time, then from the Lax equivalence theorem, 
the consistency ensures convergence.
In what follows, convergence of $\lu lF_k^n$ to $F_l(x_k,t_n)$ is proved from the global error.
It results limited by a linear, rather than exponential,
grow in time. This remarkable property is a direct consequence of the 
stochastic matrix $q_{ij}$.

\def\et{\mathcal{E}}
\begin{thm}
\label{th:conv}
Let $F_l(x,t) \in C^{2,2}(\Omega \times [0,T])$ be a solution of Eq. (\ref{ME})
with $|A_l(x)| \leq M $ and Cauchy conditions (\ref{cauchy}) for $x \in \Omega$ and $l\in\{1,\ldots,S\}$, 
where $\Omega$ is a closed set of the real axis. 
Fixed a uniform mesh on $(\Omega \times [0,T])$,
with spatial discretization of step-size $\De x=x_{k+1}-x_k$, 
and temporal discretization  $\De t=t_{n+1} - t_n$:
\begin{equation}
\De t < \left(M/\De x
+ \sup_l Q_{ll} \right)^{-1},
\label{CFL}
\end{equation}
then the global error $\normtot{\lu le_k^n}=\sum_{l=1}^S \max_k|\lu le_k^n|=\normtot{F_l(x_k,T) - \lu lF_k^n}$, of the numerical solution $\lu lF_k^n$
resulting from the scheme (\ref{upwind}),  computed at the time step $n=T/\De t$ is bounded by $\normtot{\lu le_k^n} \leq \normtot{\lu le_k^0} + T\et$,
where $\et={\mathcal O}(\De x)$, so that $\lu lF_k^n$ converges to $F_l(x_k,t_n)$ as $\De x\rightarrow 0$  
 for $(x,t) \in \Omega \times [0,T]$.
\end{thm}
\begin{pf}
We search for convergence of the numerical method by studying 
Eq. (\ref{upwind}) for the global error 
$\lu le_k^{i} = F_l(x_k,t_i) - \lu lF_k^{i}$:
\begin{equation}
\lu le_k^{i+1} = \lu le_k^{i} + \Delta t \left(-{}^l\!A_{k} 
             \frac{\lu le_{k+\nu}^{i} - \lu le_{k+\nu-1}^{i}}{\Delta x}
		+ \sum_s Q_{ls} \; \lu se_k^{i} \right) - \lu l\et_k^i\,\Delta t
\end{equation} 
where $\lu l\et_k^i$
is the local truncation error that depends on the second spatial and temporal 
derivative of the solution:
$$
\lu l\et_k^i = \frac{\De x}{2}\left(\alpha \de_t^2 \lu lF(x_k,\xi_k) -
  |\lu lA(x_k)|  \de_x^2 \lu lF(\eta_i,t_i)\right)
$$
where $\xi_k$ and $\eta_i$ are unknown points, and
 $\alpha = \De t / \De x$. Let $\lu lE^i = \max_{k}\mod{\lu le_k^{i}}$ and
$\lu l\et^i = \max_k \mod{\lu l\et_k^i}$,
taking modulus and using the triangle inequality gives:
\begin{equation}
\begin{array}{ll}
\mod{\lu le_k^{i+1}} \leq  \left( 
        \mod{ 1-\alpha|\lu lA_{k}| + \De t\, Q_{ll}} +
	 \alpha |\lu lA_{k}|\right) \lu lE^i + 
	  \De t  \sum_{s\neq l} & Q_{ls} \lu sE^i + \\
	 & \lu l\et^i\,\Delta t
\end{array}
\label{maggior}
\end{equation} 
By using the condition of Eq. (\ref{CFL}), the first modulus 
in the right side is positive:
\begin{equation}
 1-\alpha |\lu lA_{k}| + \De t\, Q_{ll}  \geq 0  \;\;\;\;\;\; \forall l,k
\label{CFL0}
\end{equation} 
then we have:
\begin{equation}
\left|\lu le_k^{i+1}\right| \leq  
       \left( 1 + \De t\, Q_{ll}  \right) \lu lE^i
	 + \De t \sum_{s\neq l} Q_{ls} \lu sE^i + \lu l\et^i\,\Delta t
\label{maggior2}
\end{equation} 
This inequality is valid for all mesh points $x_k$ and the r.h.s does not
depends on $k$, then we can write:
\begin{equation}
\lu lE^{i+1} \leq  
      \sum_{s} \left( I + \De t\, Q  \right)_{ls} \lu sE^i + \mod{\lu l\et^i}\,\Delta t
\end{equation} 
 Taking the sum over all the states $l$, we obtain 
the inequality written for the $L^1$ norm on the states:
\begin{equation}
\norm{E^{i+1}} \leq 
               \norm{( I + \De t\, Q )} \norm{E^i} 
	       + \norm{\et^i} \,\Delta t
\label{err_bound}
\end{equation}
The matrix $I + \De t\, Q $ gives us the information about the convergence
of the numerical method (\ref{upwind}).
By noting that $(I + \De t\, Q)_{ll} \geq 0$ because of (\ref{CFL0}) and 
$Q_{ls} = q_{ls}\mu_s \geq 0$ for $l\neq s$, we have:
\begin{equation}
\begin{array}{ll}
\norm{( I + \De t\, Q )} & = 
  \max_{s} \sum_{l} \vert \delta_{ls} + \De t (q_{ls} - \delta_{ls}) \mu_s)\vert \\
  & =  \max_s \left( \sum_l \delta_{ls} (1-\De t \mu_s) + 
  \De t \mu_s \sum_s q_{ls}\right) = 1
 \end{array} 
\label{norm_study} 
\end{equation}
because the fundamental property (\ref{RTE_estesa_3}) of the stochastic transition matrix $q_{ij}$.
The error of the numerical method at the time step $n$ is limited 
$$
\norm{ E^n} \leq \norm{ E^0 } + n \De t \et
$$
provided that the maximum local truncation error $\et = \max_i \norm{\et^i} $ is $\mathcal{O}(\De x)$
 for any $(x,t) \in \Omega \times [0,T]$.
\end{pf}


\subsection{Non-decreasing property}
\label{sec:non_decrease}

The solutions of Eq. (\ref{ME}) are non-decreasing function in space.
We search for the condition that the numerical solutions $\lu{l}F_{k}^{i}$
own the same property, or equivalently, that the density probability distributions,
defined as $\lu lp_k^i = (\lu lF_k^{i} - \lu lF_{k-1}^{i})/\De x$ will never becomes negative
for all $i>0$. 

\begin{thm}
Let $\lu lF_k^i$ the numerical solution obtained via the 
numerical scheme (\ref{upwind}), and  $\lu lF_k^0 - \lu lF_{k-1}^0\geq 0$.
If the CFL condition (\ref{CFL}) is satisfied,
then $\lu lF_k^i - \lu lF_{k-1}^i\geq 0$ $\forall i>0$, $\forall l, \forall k$.
\end{thm}

\begin{pf}
By subtracting the solution at meshes point $k$ and $k+1$, we find: 
\begin{equation}
\begin{array}{l}
\lu{l} F_k^{i+1} - \lu{l} F_{k-1}^{i+1} = \lu{l}F_{k}^{i} - \lu{l} F_{k-1}^{i} +
             \De t \left( 
	     - \frac{\lu{l} A_k}{\De x} \lu{l} F_{k+\nu}^{i}
	     + \frac{\lu{l} A_k}{\De x} \lu{l} F_{k+\nu-1}^{i} \right. \\
	     \left. 
	     + \frac{\lu{l} A_{k-1}}{\De x} \lu{l} F_{k+\nu-1}^{i} 
	     - \frac{\lu{l} A_{k-1}}{\De x} \lu{l} F_{k+\nu-2}^{i} 
	     + \sum_s Q_{ls} \left( \lu{s}F_k^{i} - \lu{s}F_{k-1}^i 
	     \right) \right) 	    
\end{array}
\end{equation}

Let $\De x\,\lu lp_k^{i+1} = \lu lF_k^{i+1} - \lu lF_{k-1}^{i+1}$, we have:

\begin{equation}
\lu{l}p^{i+1}_{k} = \lu{l}p^{i}_{k} + \De t \left( 
                 \frac{|\lu{l}A_{k+\nu-1}|}{\De x} \lu{l}p_{k+2\nu-1}^{i} +
		 \sum_s \left( Q_{ls} - \frac{|\lu{l}A_{k-\nu}|}{\De x}
		 \delta_{ls}
		 \right)  \lu{s}p^i_k \right) 
\end{equation} 
where $Q_{ls}=(q_{ls}-\delta_{ls})\mu_s$, that is:

\begin{equation}
\begin{array}{ll}
\lu{l}p^{i+1}_{k} =& \left( 1 - \De t \left( 
                 \frac{|\lu{l}A_{k-\nu}|}{\De x} -Q_{ll} \right) \right)
		 \lu{l}p_{k}^{i} +
		  \De t \sum_{s\neq l} \left( Q_{ls} \lu{s}p^i_k \right) + \\
		 & \frac{\De t}{\De x}|\lu{l}A_{k+\nu-1}| \lu{l}p^i_{k+2\nu-1} 
\end{array}
\end{equation}
We have to ensure that $\lu{l}p^{i+1}_{k}$ will never becomes negative. 
Suppose that exists $i$ such that $\lu{l}p^{i}_{k}\geq 0$ for all $k,l$, then
$$
\lu lp_k^{i+1} \geq \left( 1 - \De t \left(\mod{\lu lA_{k-\nu}}/\De x -Q_{ll} \right) 
\right) \lu lp_{k}^{i}
$$
because of $Q_{ls}\geq 0$ for $l\neq s$. Using the CFL restriction for $\De t$, we find
$\lu lp_k^{i+1} \geq \lu lp_k^{i} \geq 0$. Now, being $\lu lp_k^0\geq 0$, by induction we find 
the thesis. 
\end{pf}
We observe the CFL condition is a sufficient condition for both convergence and 
non-decreasing of the numerical solution $\lu{l}F^{i}_{k}$ (or 
non negativity of $\lu{l}p^{i}_{k}$).

\subsection{Smearing}
In the theorem \ref{th:conv} we have restricted our analysis to 
solutions having both second spatial and temporal derivative,
but in general local discontinuities can be considered because of 
the upwind method tends to regularize the solution.
This is explained in terms of \textit{modified equation}
\cite{ran} modeling the numerical method.
By substituting a function $F(x,t)$ into Eq. (\ref{upwind}),
expanding in Taylor's series and preserving the second order terms, in 
matrix notation we get:
\begin{equation}
\de_t F + A(x)\de_x F + \frac{\De t}{2} {\de_t^2 F} 
           - \frac{A(x) \De x}{2} \de_x^2 F - Q F \simeq 0 .
\label{second_order}
\end{equation}
Taking the derivatives by  $\de_x$ and $\de_t$ of this equation, we can write:
$$
\de_t^2 F \simeq A(\de_x A-Q)\de_x F + Q \de_t F + A^2 \de_x^2 F
$$
that inserted into Eq. (\ref{second_order}) give us:
\begin{equation}
\begin{array}{ll}
\left( I - \frac{\De t}{2} Q\right) \de_t F + 
      A(x) \left[I+\frac{\De t}{2}(\de_x A - Q)\right] & \de_x F -Q F =  \\
      &  \frac{A(x) \De x}{2} \left( I - \frac{A(x)\De t}{\De x} \right) \de_x^2 F
\end{array}
\label{modified}      
\end{equation}
When $\Delta t \rightarrow 0$ we recover Eq. (\ref{ME}) showing 
the consistency of this modified equation.
Coefficients of $\de_t$ and $\de_{xx}$ are diagonal matrices and Eq. (\ref{modified}) 
owns, loosely speaking, the diffusive properties of a parabolic equations.
The upwind method (\ref{upwind}) is of the second order with respect to
Eq. (\ref{modified}), so that the numerical solution $\lu lF_k^i $ is closer to this
modified equation than Eq. (\ref{ME}), showing a diffusive character.
%


\section{Numerical tests}
\label{numtest}

We perform numerical tests by assigning some values to the parameters model.
Then we solve the Cauchy problem for the 
evolutory equation (\ref{ME}), using the upwind method (\ref{upwind}), 
and plot the total density probability $p_k^i = \sum_{l=1}^S \lu lp_k^i$. 

\begin{figure}
   \centering
   \includegraphics[scale=.9]{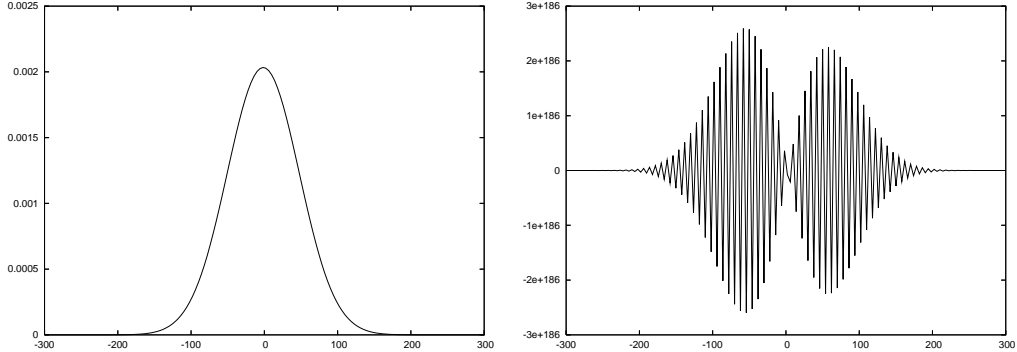}
   \caption{On the left: convergent solution $p(x_k,T)$ of Eq. (\ref{ME}) when
   the CFL condition (\ref{CFL}) is satisfied. On the right: non-convergent
   solution when that is not satisfied. }
   \label{fig:check_CFL}
\end{figure}

\begin{figure}[!t]%
\includegraphics[scale=.54]{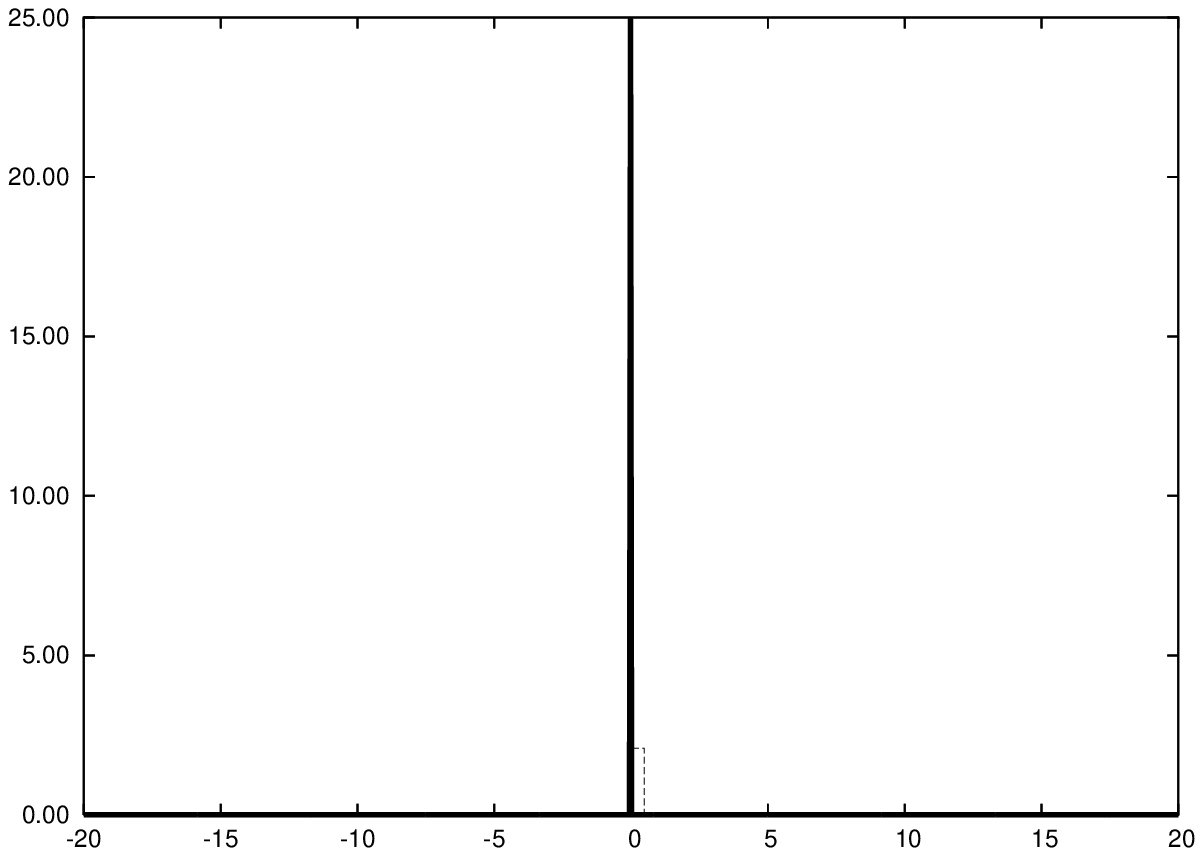}%
\includegraphics[scale=.54]{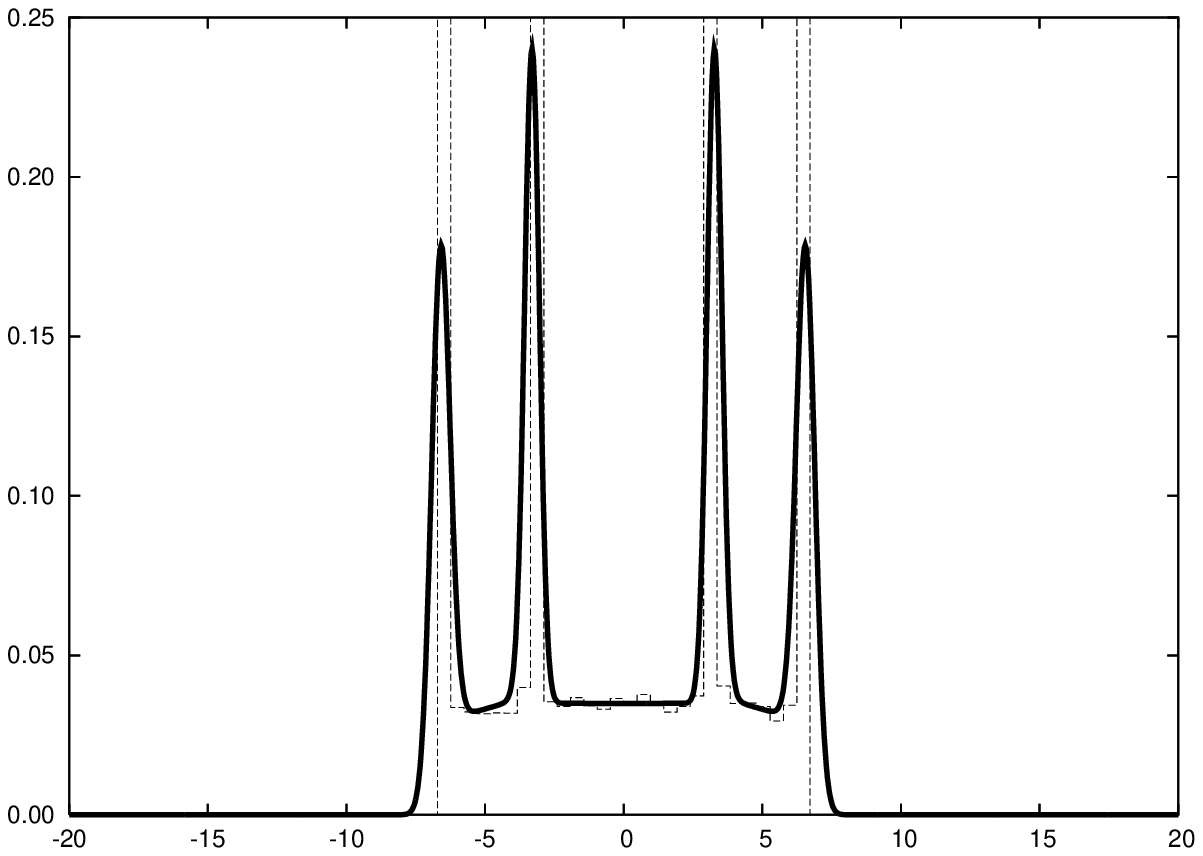}
\includegraphics[scale=.54]{./figure/broad_art.epsi}%
\includegraphics[scale=.54]{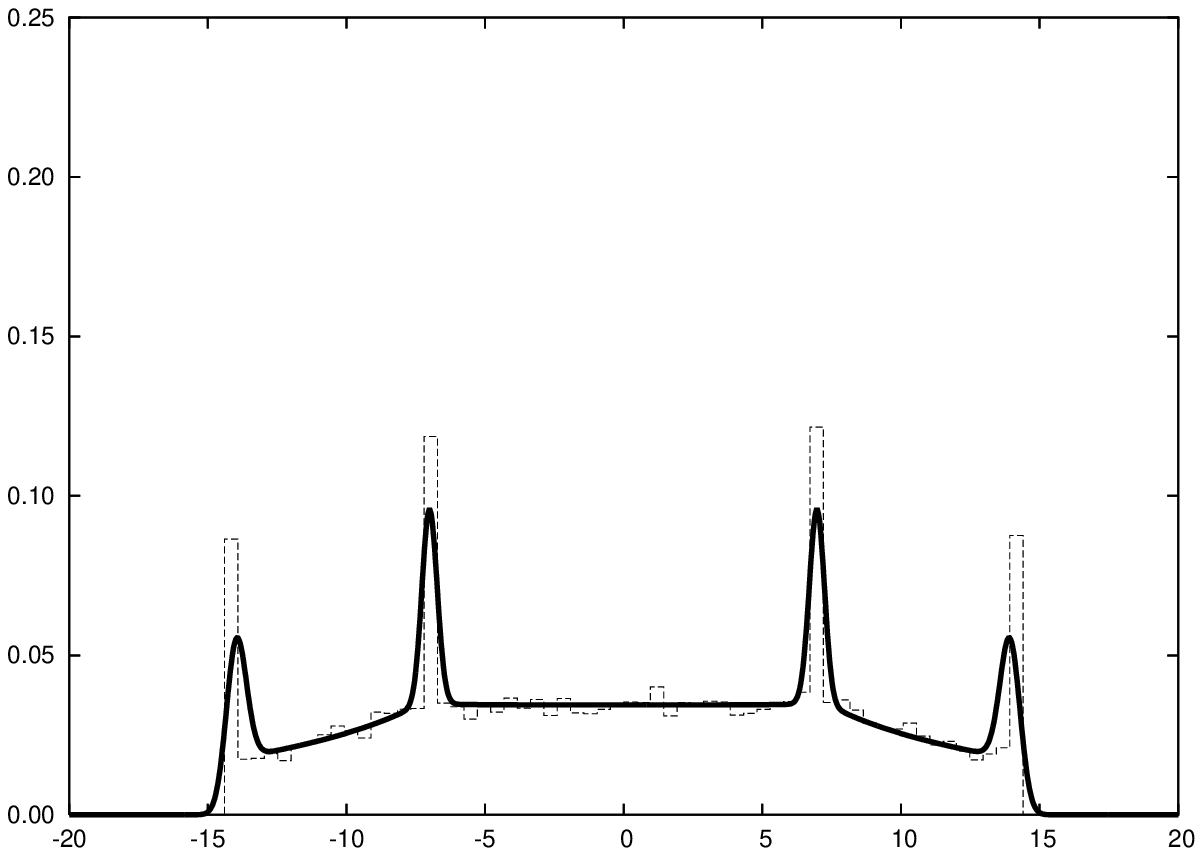}
\includegraphics[scale=.54]{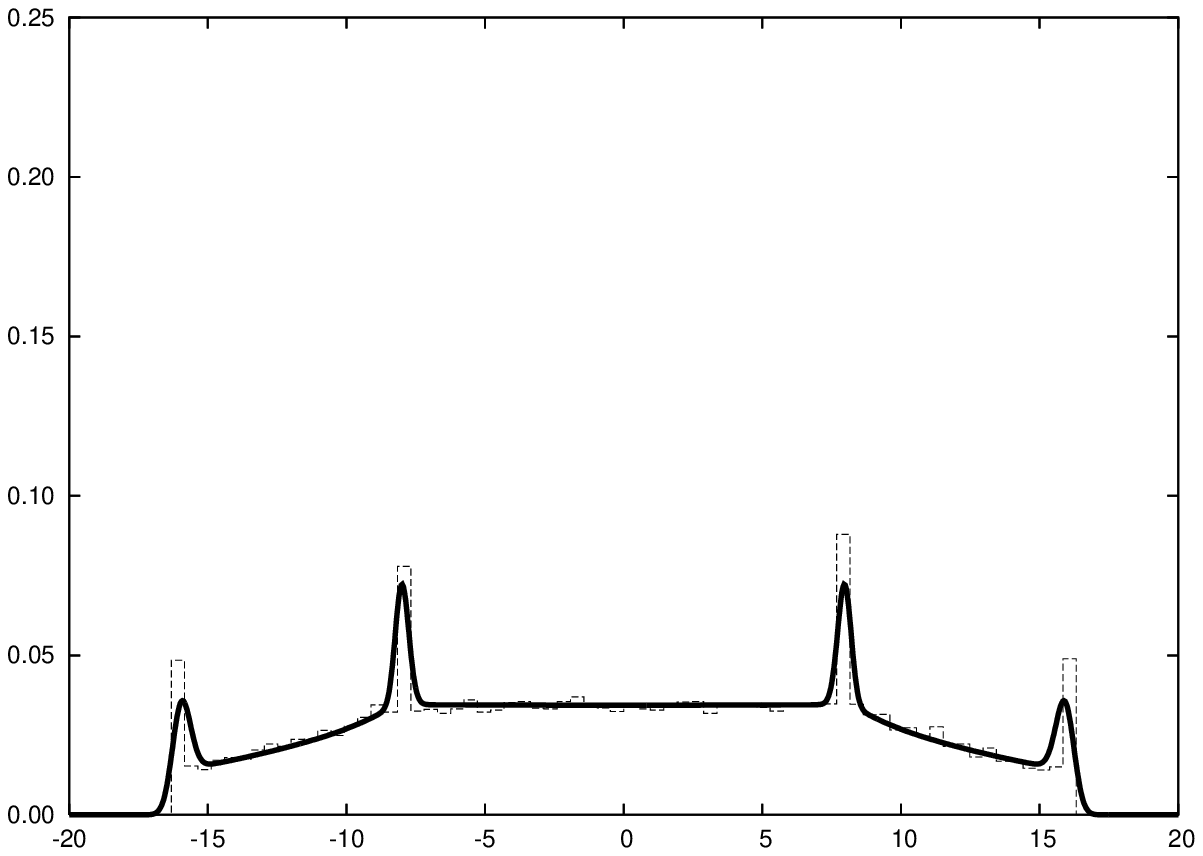}%
\includegraphics[scale=.54]{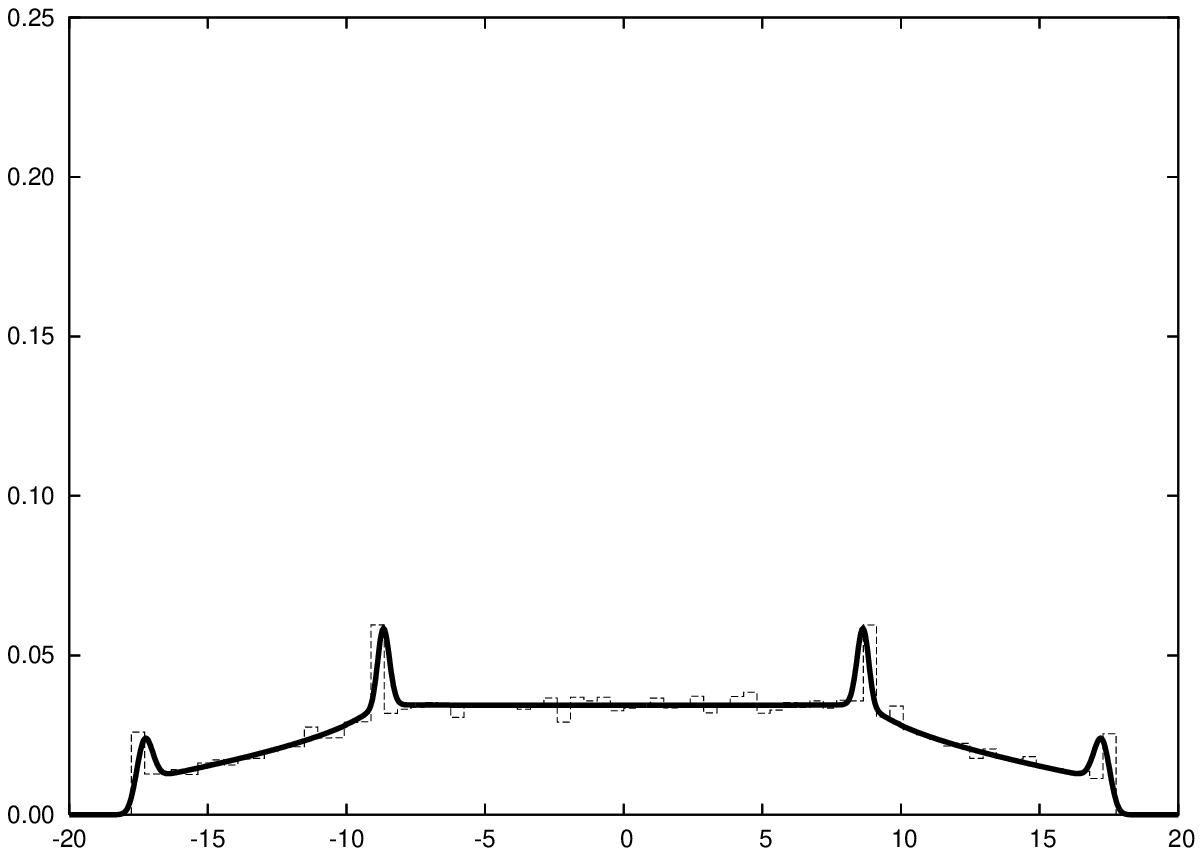}
\caption{Temporal evolution of the density probability distribution functions $p(x,t)$ at fixed time. $S=4$, $\De t=0.009$, $\De x = 0.04004$, $\mu_j =0.2$,
$\gamma_j=0.1$. From top left to bottom right: $t=\{0,4,8,12,16,20\}$.
Monte Carlo's histograms in dashed lines.} 
\label{temporal_evolution}
\end{figure}

In order to test the validity of the numerical scheme versus
the CFL condition, we study a system of a fluctuation-dissipation process having four relaxation states:
$A_s(x) = -\gamma_s x + W_s$  with $s=\{1,2,3,4\}$, where 
the switching times happen according to the Poisson processes
of Eq. (\ref{RTE_estesa_2}).
This model is solved for the initial condition (Riemann problem) 
$F_{0s}(x) = H(x)/4$ (or $p_{0s}(x) = \delta(x)/4$), up to  
time $T$, where $H(x)$ is the Heaviside function.
As discussed in Sect. \ref{sec:extend} this kind of system the density distribution always reaches an equilibrium, i.e. 
$\lim_{t\rightarrow\infty} p(x,t) = p_{eq}(x)$ with
$x \in \Omega \equiv \left[ \min_{s}\left\lbrace  \frac{W_s}{\gamma_s}\right\rbrace; 
\max_{s} \left\lbrace \frac{W_s}{\gamma_s}\right\rbrace \right]$.
The discontinuous initial condition does not fulfill the hypothesis for convergence, but, as above mentioned,  it is subject to smearing 
that tends to reguralize discontinuities 
(see the third picture of Fig. (\ref{temporal_evolution})).
In Fig. (\ref{fig:check_CFL}) we show the total density distribution function 
$p(x_k,T)$ when the temporal step $\De t$ meet or not the CFL condition.
We use the following parameters:
$\mu_s = 4$, $q_{ls}=0.25$, $\gamma_s = 10^{-3}$, 
$W_1=1, W_2=-1, W_3=2, W_4=-2$. By using $\De x = 4.004$, we 
get from the CFL condition of Eq. (\ref{CFL}) that 
$\De t < \De t_{max} = 0.250063$.
On the left side of Fig. (\ref{fig:check_CFL}), we used $\De t = 0.225056 < \De t_{max}$
and $T= 545.5$. On the right side: $\De t = 0.5 > \De t_{max}$ and the 
solution is not convergent.
In Fig. (\ref{temporal_evolution}) we plot the total probability $p(x_k,t_i)$ 
for fixed times to see how the density distribution evolves. 
For these plots $\mu_j =0.2$,
$\gamma_j=0.1$, the remaining parameters are same to previous test.
These results are compared to histograms
obtained by Monte Carlo simulation of the dynamical equations
(\ref{RTE_estesa_1}),(\ref{RTE_estesa_2}) and (\ref{RTE_estesa_3}).
Histograms are generated by collecting the final values $X(t)$
of the sample paths.

\section{Summary and conclusions}

In this communication we have introduced the continuous piecewise deterministic  Markov processes, 
as an extention of random telegraph processes, and the related Liouville-Master Equation 
for the distribution functions.
The upwind method has been successfully applied for solving this equation.
We have found a Courant-Friederichs-Lewy condition ensuring both convergence  and
non-decreasing property of the numerical distribution function. 
The numerical solution fits with that we obtained by performing Monte Carlo
simulation of the process.
These results are an introductory work for this general class of processes.
We trust that the piecewise deterministic processes have a potential, as
models for non-deterministic systems, large range of application 
as confirmed by the recent literature involved on the subject,
thus justifying the development of numerical methods for Eq. (\ref{ME}).
Many aspects should be further investigated, such as: the asymptotic stability for stiff parameters, 
the behaviour for discontinuous distributions, the development of higher order numerical methods,
the extention of the model to non-Markovian and many dimensional processes.


\begin{thebibliography}{00}




\bibitem{DAV} M.H.A. Davis,\textit{Markov Models and Optimization}, Monograph on 
Statistics and Applied Probability 49, Chapman \& Hall/CRC 1993.

\bibitem{gar} C.W. Gardiner,
\textit{Handbook of Stochastic Methods}, Springer-Verlag Berlin Heidelberg New York 
2nd Ed. 1985.

\bibitem{mor:may} K.W. Morton \& D.F. Mayers,
 \textit{Numerical Solution of Partial Differential Equations}, Cambridge University Press 1994.

\bibitem{ran} Randall J. LeVeque, \textit{Numerical Methods for conservation Laws},
Birkh\"auser 1992.
 

\bibitem{mario} M. Annunziato, ``Non-gaussian equilibrium distribution arising from the Langevin equation," Phys. Rev. E, {\bf 65} 21113 (2002).

\bibitem{bol:gri:wes} M. Bologna, P. Grigolini, B.J. West, ``Strange kinetics:
conflict between density and trajectory description,'' Chem. Phys. {\bf 284} 115 (2002).

\bibitem{bud:cac} A. Budini, M.O C\'aceres, ``Functional characterization of generalized
Langevin equations,'' J. Phys. A: Math. Gen. {\bf 37} 5959 (2004).

\bibitem{cac}
M.O C\'aceres, ``Computing a non-Maxwellian velocity distribution
from first principles,'' Phys. Rev E {\bf 67} 016102 (2003)
 
\bibitem{dav}
M.H.A. Davis, ``Piecewise-Deterministic Markov Processes:
A General Class of Non-Diffusion Stochastic Models,'' J. of the Royal
Stat. Soc. Series B {\bf 46} 353 (1984).

 

\bibitem{fil:hon}
R. Filliger, M.O. Hongler, ``Supersymmetry in random two-velocity 
processes,'' Physica A {\bf 332} 141 (2004).


\bibitem{oth:dun:alt} H. G. Othmer, S. R. Dunbar, and W. Alt, ``Models of dispersal in biological systems,'' J.
Math. Biol. {\bf 26} 263 (1988).


\bibitem{hor} W. Horsthemke, 
``Spatial instabilities in reaction random walks with 
direction-independent kinetics,'' Phys. Rev E {\bf 60} 2651 (1999).

\bibitem{jak:ren} E. Jakeman, E. Renshaw, ``Correlated random-walk model for scattering," J. Opt. Soc. Am. A {\bf 4} 1206 (1987).


\bibitem{jak:rid} E. Jakeman,  K.D. Ridley,
 ``The statistics of a filtered telegraph signal," J. Phys. A, {\bf 32} 8803 (1999).

\bibitem{kit} K. Kitahara \textit{et al.}, ``Phase Diagrams of Noise Induced Transitions,'' Prog. Theor.
Phys. {\bf 64} 1233 (1980).


\bibitem{mas:por} J. Masoliver, J. M. Porr\`a, ``Harmonic oscillators driven by colored noise:
Crossovers, resonances, and spectra," Phys. Rev. E {\bf 48} 4309 (1993).

\bibitem{mor} A. Morita, ``Free Brownian motion of a particle driven
by a dichotomous random force,'' Phys. Rev. A {\bf 41} 754 (1990).

\bibitem{paw:ric} 
R.F. Pawula, O. Rice,``On Filtered Binary Processes," IEEE Trans. Inf. Th. Vol. IT-32, {\bf 63} (1986).

\bibitem{bro:han} C. Van den Broeck, P. H\"anggi, ``Activation rates for nonlinear stochastic flows driven by non-Gaussian noise,'' Phys. Rev. A {\bf 30} 2730 (1984).

\bibitem{won} W.M. Wonham,  J. Electron. Control {\bf 6} 376 (1959).

\end{thebibliography}
\end{document}